\documentclass[numberwithin=section]{papier}

\title{Atomic Toposes with Co-Well-Founded Categories of Atoms}
\author{Jérémie Marquès}
\date{\today}

\newcommand{\at}{\mathrm{at}}

\bibliography{references.bib}

\begin{document}
\maketitle

\textbf{Disclaimer:} This is an earlier version of an article published in open access in \emph{Theory and Applications of Categories} \cite{MarAtomicToposesCoWellFounded2025} (\href{http://www.tac.mta.ca/tac/volumes/44/15/44-15abs.html}{link}). I leave this version here because I removed the discussion of Fraïssé limits in §~\ref{subsec:pBpresSheaves}, but it may still have some interest. It is currently a bit sketchy and too complex for a result whose interest is unclear. In \cite{MarAtomicToposesCoWellFounded2025}, Proposition~\ref{prop:all-atoms} has been improved with a complete description of the category of atoms; this answers the first question in §~\ref{sec:final-rmk} of the current version. An example of an atomic topos which is not locally finitely presentable has also been added.

\begin{abstract}
	The atoms of the Schanuel topos can be described as the pairs $(n,G)$ where $n$ is a finite set and $G$ is a subgroup of $\Aut(n)$. We give a general criterion on an atomic site ensuring that the atoms of the topos of sheaves on that site can be described in a similar fashion. We deduce that these toposes are locally finitely presentable. By applying this to the Malitz--Gregory atomic topos, we obtain a counter-example to the conjecture that every locally finitely presentable topos has enough points. We also work out a combinatorial property satisfied exactly when the sheaves for the atomic topology are the pullback-preserving functors. In this case, the category of atoms is particularly simple to describe.
\end{abstract}

\showcontents
\vspace{1em}

The original motivation for this article was the question posed in \cite{dilibertiTopoiEnoughPoints2024} about whether every locally finitely presentable topos has enough points. We show here that the pointless atomic topos given in \cite[Sec.~5]{makkaiFullContinuousEmbeddings1982} is a counter-example. In order to do so, we give a combinatorial criterion on an atomic site ensuring that the atoms of its topos of sheaves are not ``too complex,'' and as a consequence that this topos is locally finitely presentable. A crucial property is that these categories of atoms are \emph{co-well-founded} in the sense that every chain stabilizes.

As a motivating example, we consider the Schanuel topos of nominal sets. It can be defined in many ways; for a complete introduction, we refer to the book \cite{NominalSetsNamesPit2013}. Informally, an object of this topos, a nominal set, can be thought of as a set of \emph{terms} in which a finite number of \emph{variables} may appear. The variables of a term can be renamed to yield what is called in computer science an $α$-equivalent term.


In order to give an idea of what nominal sets look like in general, we examine two examples.
\begin{enumerate}
	✦ The nominal set of \emph{ordered pairs} has terms $(a,b)$ where $a$ and $b$ are variables.
	✦ The nominal set of \emph{unordered pairs} has terms $\set{a,b}$ where $a$ and $b$ are variables.
\end{enumerate}
The difference between these two nominal sets is that exchanging the variables of an unordered pair does nothing, since $\set{a,b} = \set{b,a}$ but $(a,b) ≠ (b,a)$ if $a≠b$. These two examples actually illustrate everything that can happen in a nominal set. More precisely, each nominal set can be decomposed as a potentially infinite union of basic building blocks, its \emph{atoms}. An atom of a nominal set is a term modulo renaming of its variables. This means that we can think of an atom as a term $t(x_1,…,x_n)$ where the $x_i$ are pairwise distinct variables. The number of variables is the \emph{support} of the atom. Another important information is the group $G ⊆ \gSym_n$ of permutations of the variables $x_i$ that leave the term unchanged. The atom is in a sense completely determined by the pair $(n,G)$. In categorical terms, an atom of a nominal set is a minimal non-empty sub-object, and two atoms are isomorphic if and only if they have the same invariant $(n,G)$.

Let us determine for instance the atoms of the nominal set of unordered pairs. Modulo renaming, there are two unordered pairs: $\set{a,b}$ with $a≠b$ and $\set{a}$. The invariants of these atoms are respectively $(2,\gSym_2)$ and $(1,\gSym_1)$.

The Schanuel topos can also be presented as the category of sheaves for the atomic topology on the opposite of the category $\cFinSetInj$ of finite sets and injections. There is a direct connection with the explicit description of the atoms given above: each atom is a formal quotient of an object $n$ of $\cFinSetInj^\op$ by a group of automorphisms of $n$. In this paper, we will give a condition on a category $\cC$ so that this phenomenon generalizes to the topos of sheaves for the atomic topology on $\cC^\op$.

In §~\ref{sec:well-founded}, we show that if the category of atoms of an atomic topos is co-well-founded, then it is locally finitely presentable. We state in §~\ref{sec:general-crit} our criterion to show that an atomic topos has a co-well-founded category of atoms. We focus in §~\ref{subsec:pBpresSheaves} on a case where the category of atoms is particularly simple to describe. Finally, we apply the criterion in §~\ref{sec:exmp-malitz-gregory} to show that a classical example of pointless atomic topos, the Malitz--Gregory topos, is locally finitely presentable. This gives a counter-example to the conjecture that every locally finitely presentable topos has enough points formulated in \cite{dilibertiTopoiEnoughPoints2024}.

\paragraph{Acknowledgments} I would like to thank Morgan Rogers for the discussions about his recent paper \cite{dilibertiTopoiEnoughPoints2024} with Ivan Di~Liberti, from which this document originated. I am also grateful to my post-doc advisor Sam van~Gool for his many comments and suggestions. Finally, my intuitions about nominal sets and atomic toposes greatly benefited from discussions with Victor Iwaniack and Daniela Petrişan. The research reported here has been supported financially by the European Research Council (ERC) under the European Union's Horizon 2020 research and innovation program, grant agreement \#670624.

\paragraph{Notations} In this paper, topos means sheaf topos. The symbol $J_\at$ denotes the atomic topology on a category left implicit in the notation. The composition of two arrows $f : A→B$ and $g : B→C$ is denoted by $fg$ or $g∘f$.

\section{Co-well-founded categories of atoms and finite presentability}
\label{sec:well-founded}

A \emph{chain} in a category is a diagram indexed by a well-ordered poset. We say that a chain $(x_i)_{i∈I}$ \emph{stabilizes} at $i∈I$ if the morphisms $x_i→x_j$ are isomorphisms for all $j > i$. If a chain stabilizes at $i$, then it stabilizes at every $j > i$ and the morphism $x_j → \colim_k x_k$ is an isomorphism for all $j > i$. An \emph{$ω$-chain} is a chain indexed by $ω$. We say that a category $\cC$ is \emph{co-well-founded} if every chain stabilizes. Equivalently, $\cC$ is co-well-founded if every $ω$-chain stabilizes: If there is a chain that doesn't stabilize, pick an element $a_0$ in it. Then there is a morphism $a_0 → a_1$ in the chain which isn't an isomorphism. Repeating this process produces an $ω$-chain that doesn't stabilize. We say that $\cC$ is \emph{well-founded} if $\cC^\op$ is co-well-founded.

\begin{rmq}{}{}
	Another possible name for a co-well-founded category could be ``Noetherian category.'' However, this usually means that every chain of sub-objects of a any given object stabilizes.
\end{rmq}

An atom of an object in an atomic topos $\Es$ is a minimal non-empty sub-object. Let $π_0 : \Es→\cSet$ be the functor taking an object to its set of atoms. We say that $X$ is \emph{finite} if $π_0(X)$, or equivalently $\Sub_\Es(X)$, is finite.

\begin{prop}{}{stab-lfp}
	If the category of atoms of an atomic topos is co-well-founded, then the topos is locally finitely presentable. In this case, the finitely presentable objects are exactly the finite objects.
\end{prop}

\begin{proof*}{}
	In an atomic topos, the class of atoms is essentially small, see \cite[p.~690]{johnstoneSketchesElephantTopos2002}. Moreover, any object is the filtered colimit of its finite sub-objects. It remains only to show that, under the hypothesis of the property, for any atom $a$, the functor $\Hom(a,—)$ preserves filtered colimits. By \cite[Cor.~1.7, p.~15]{adamekLocallyPresentableAccessible1994}, it suffices to show that $\Hom(a,—)$ preserves colimits of chains. Let $(X_i)_{i∈I}$ be a chain in $\Es$ and let $X = \colim_i X_i$. Let $a → X$ be a morphism. The union of the images of $X_i → X$ is $X$, hence there is an $i$ and an atom $b ∈ X_i$ whose image coincides with the image of $a$ in $X$. The sequence of images of $b ∈ X_i$ in $X_j$ for $j > i$ stabilizes since the category of atoms is co-well-founded. Hence, $a→X$ factors through $X_j → X$ for some $j > i$. Now, consider two different morphisms $a ⇉ X_i$ such that the composites $a ⇉ X_i → X$ are equal. There is an index $j > i$ such that the two morphisms $a ⇉ X_j$ have the same \emph{image} $b ⊆ X_j$, since $π_0$ is cocontinuous. Since the sequence of atoms obtained as the images of $b$ in $X_k$ for $k > j$ stabilizes, there is some $k > j$ such that the two morphisms $a⇉X_k$ are equal. This shows that $\Hom(a,\colim_i X_i) ≅ \colim_i \Hom(a,X_i)$.
\end{proof*}

\begin{rmq}{}{}
	We will sketch a different proof of Proposition~\ref{prop:stab-lfp}. Given a category $\cC$, let $⨆[\cC]$ be its free cocompletion under small coproducts and let $⊔[\cC]$ be its free cocompletion under finite coproducts. We have an equivalence $⨆[\cInd(\cC)] ≃ \cInd(⊔[\cC])$. If the category $\cA$ of atoms is co-well-founded, then $\cA ≃ \cInd(\cA)$ and we have $\Sh(\cA,J_\at) ≃ ⨆[\cA] ≃ ⨆[\cInd(\cA)] ≃ \cInd(⊔[\cA])$.
	
	This could be put in comparison with \cite{borceuxLeftExactPresheaves1999}, where a criterion involving well-foundedness on a pretopos $\Es$ is given so that $\cInd(\Es)$ is the topos of sheaves on $\Es$ with the coherent topology (in general, $\cInd(\Es)$ is only included in this topos).
\end{rmq}

\begin{exmp}{}{}
	In general, every finitely presentable object of an atomic topos is an atom, but the converse isn't true. There are for instance locally finitely presentable atomic toposes whose categories of atoms are not co-well-founded, or equivalently which possess non finitely presentable atoms. To build an example, let $F_2$ be the free group on two generators $x$ and $y$. By the Nielsen--Schreier Theorem, every subgroup of $F_2$ is free. Let $G$ be the kernel of the homomorphism $F_2→ℤ$ sending $x$ and $y$ to $1$. It is freely generated by all the elements of the form $x^ny^{-n}$ for $n∈ℤ$. In the topos of $F_2$-sets, the atom $F_2∕G$ is not finitely presented since it is the colimit of $F_2∕\genrel{x^ny^{-n}}{-N < n < N}$ as $N$ goes to infinity.
\end{exmp}

\section{A well-foundedness criterion}
\label{sec:general-crit}

In this section, we give a criterion on a small category $\cC$ to ensure that the atoms of $\Sh(\cC^\op,J_\at)$ form a co-well-founded category. It is composed of the conditions \ref{cond:1}, \ref{cond:2}, \ref{cond:3}, \ref{cond:4} below. We express here \ref{cond:1} and \ref{cond:2} in an abstract form, but we will give equivalent combinatorial formulations in §~\ref{subsec:combinatorial-content}.
\begin{itemize}
	✦[\optionaldesc{(C1)}{cond:1}] $\cC^\op → \Sh(\cC^\op,J_{\at})$ is fully faithful and sends each object of $\cC^\op$ to an atom.
	✦[\optionaldesc{(C2)}{cond:2}] $\cC^\op$ has pushouts and $\cC^\op → \Sh(\cC^\op,J_{\at})$ preserves them.
	✦[\optionaldesc{(C3)}{cond:3}] $\cC$ is well-founded.
	✦[\optionaldesc{(C4)}{cond:4}] The automorphism groups of the objects of $\cC$ are Noetherian.
\end{itemize}
It seems there that a formulation in terms of $\cC^\op$ instead of $\cC$ would be more natural. However, we stick with $\cC$ because it can be described in the examples as a category of finitary models of some theory, while $\cC^\op$ is harder to work with directly. In terms of the syntax-semantics duality, $\cC$ is on the semantical side while $\cC^\op$ is on the syntactical side. 

The plan is as follows. Condition \ref{cond:1} ensures that we can describe a category of ``representable'' atoms equivalent to $\cC^\op$. Condition \ref{cond:2} ensures that we can compute pushouts of these atoms. Condition \ref{cond:3} will be used to compute coequalizers of morphisms between representable atoms by an iteration of pushouts followed by a quotient by an automorphism. From this, we deduce a description of all the atoms of $\Sh(\cC^\op,J_\at)$. Finally, we use \ref{cond:4} to show that this atomic topos has a co-well-founded category of atoms.

\subsection[The combinatorial content of (C1) and (C2)]{The combinatorial content of \ref{cond:1} and \ref{cond:2}}
\label{subsec:combinatorial-content}

We say that a category $\cC$ has \emph{amalgamation} or satisfies the \emph{left Ore condition} if every span admits a cocone.

As can be found in \cite[Sec.~7(3)]{AtomicToposesBarDia1980} and \cite[C2.1.12(c)]{johnstoneSketchesElephantTopos2002}, \ref{cond:1} is equivalent to the conjunction of
\begin{enumerate}
	✦ $\cC$ has amalgamation, and
	✦ every morphism in $\cC$ is the equalizer of all the pairs of maps that it equalizes.
\end{enumerate}
When $\cC$ has these properties, $(\cC^\op,J_\at)$ is called a \emph{standard} atomic site in \cite{AtomicToposesBarDia1980}. In this case, we call the atoms in the image of $\cC^\op → \Sh(\cC^\op,J_\at)$ the \emph{representable atoms}.

In order to give the combinatorial content of \ref{cond:2} in Lemma~\ref{lem:homogen-combinatorial}, we introduce another condition on $\cC$:
\begin{itemize}
	✦[\optionaldesc{(C2')}{cond:2-}] $\cC$ has pullbacks and for every pullback
	\begin{equation}\label{diag:pb}\begin{tikzcd}
			X∩Y \ar[r] \ar[d] & X \ar[d,"f"] \\
			Y \ar[r,"g"] & Z
	\end{tikzcd}\end{equation}
	in $\cC$, for every pair of parallel morphisms $u,v : Z⇉A$ coinciding on $X∩Y$, there is a morphism $w : A→A'$ and a sequence of morphisms $k_0 = uw, k_1, k_2, …, k_n = vw : Z→A'$ such that any two consecutive morphisms in this sequence coincide either on $X$ or on $Y$.
\end{itemize}

If $\cC$ satisfies \ref{cond:1}, then every arrow is a monomorphism. This justifies the notation $X∩Y$, but in general it is just a pullback.

\begin{lem}{}{homogen-combinatorial}
	Let $\cC$ be a small category. Then \ref{cond:2} and \ref{cond:2-} are equivalent.
\end{lem}

\begin{proof*}{}
	The canonical functor $y : \cC^\op → \Sh(\cC^\op,J_\at)$ preserves pushouts if and only if every sheaf $F : \cC→\cSet$ preserves pullbacks, using that $\Hom(y(n), F) ≅ F(n)$. Suppose that $\cC$ satisfies \ref{cond:2-}. Let $F$ be a sheaf for the atomic topology. We show that it preserves the pullback \eqref{diag:pb}. We have $F(X∩Y) ⊆ F(X)∩F(Y)$ and we must show the reverse inclusion. Let $z ∈ F(X)∩F(Y)$. We show that $z$ satisfies the sheaf covering condition for $X∩Y→Z$. Let $u,v : Z⇉A$ be a parallel pair of arrows coinciding on $X∩Y$. We must show that $u(z) = v(z)$. Let $w : A→A'$ and $k_0 = uw, k_1, k_2, …, k_n = vw$ a sequence of morphisms as in \ref{cond:2-}. Then $k_0(z) = k_1(z) = ⋯ = k_n(z)$. This shows that $z$ descends to $X∩Y$, and $\cC$ satisfies \ref{cond:2}.
	
	Suppose now that $\cC$ satisfies \ref{cond:2} and consider the square \eqref{diag:pb}. Let $P : \cC→\cSet$ be the pushout of $\Hom(X,—)$ and $\Hom(Y,—)$ along $\Hom(Z,—)$ computed in $[\cC,\cSet]$. Let $s : [\cC,\cSet] → \Sh(\cC^\op,J_\at)$ be the sheafification functor. Then \ref{cond:2} says that $s(P) → s(\Hom(X∩Y,—))$ is an isomorphism.
	\begin{equation}\label{diag:sheafify-pushout}\begin{tikzcd}
		P \ar[r] \ar[d] & \Hom(X∩Y,—) \ar[d] \\
		s(P) \ar[r] & s(\Hom(X∩Y,—))
	\end{tikzcd}\end{equation}
	Let $u, v : Z⇉A$ be two arrows that coincide on $X∩Y$. These arrows represent two elements $[u], [v]$ of $P(A)$ which are sent to the same element of $\Hom(X∩Y,A)$. Using the commutativity of \eqref{diag:sheafify-pushout} and the fact that $s(P) → s(\Hom(X∩Y,—))$ is an isomorphism, we obtain that $[u]$ and $[v]$ are sent to the same element of $s(P)(A)$. This means that $[u]$ and $[v]$ are locally equal, i.e. that there is a morphism $w : A→A'$ such that $[uw] = [vw]$ in $P(A')$. The definition of $P(A')$ as a pushout gives the sequence of morphisms in \ref{cond:2-}.
\end{proof*}

As we saw in the proof, \ref{cond:2} implies that every sheaf is a pullback-preserving presheaf. We will come back to this in in §~\ref{subsec:pBpresSheaves}, where we consider the converse implication.

\subsection{The atoms of $\Sh(\cC^\op,J_\at)$}

Given an object $A ∈ \Sh(\cC^\op,J_\at)$ and an automorphism $σ$ of $A$, we denote by $A∕σ$ the coequalizer of $σ$ and the identity. If $G ⊆ \Aut(m)$ is a subgroup, $A∕G$ denotes the common coequalizer of all the arrows in $G$.

\begin{lem}{}{coeq}
	Suppose that $\cC$ satisfies \ref{cond:1}, \ref{cond:2} and \ref{cond:3}. Let $n⇉m$ be a pair of morphisms between representable atoms in $\Sh(\cC^\op,J_\at)$. Then their coequalizer is of the form $m→m'→m'∕σ$ where $m'$ is a representable atom and $σ$ is an automorphism of $m'$.
\end{lem}

\begin{proof*}{}
	Let $α_0,β_0 : m_0⇉m_1$ be a pair of morphisms between representable atoms. Form the pushout below.
	\[\begin{tikzcd}
		m_0 \ar[r,"α_0"] \ar[d,"β_0"] & m_1 \ar[d,"α_1"]\\
		m_1 \ar[r,"β_1"] & m_2
	\end{tikzcd}\]
	The coequalizer of $α_0$ and $β_0$ is canonically isomorphic to the coequalizer of $α_1$ and $β_1$. We iterate the process and we obtain a sequence of representable atoms
	\[\begin{tikzcd}
		m_0 \ar[r,"α_0",shift left=0.4em] \ar[r,"β_0"',shift right=0.4em] & m_1 \ar[r,"α_1",shift left=0.4em] \ar[r,"β_1"',shift right=0.4em] & m_2 \ar[r,"α_2",shift left=0.4em] \ar[r,"β_2"',shift right=0.4em] & ⋯
	\end{tikzcd}\]
	Since every chain of representable atoms stabilizes by \ref{cond:3}, there is some $i$ such that $α_i$ and $β_i$ are isomorphisms. The coequalizer of $α_i$ and $β_i$ is $m_i∕(α_iβ_i^{-1})$, hence the coequalizer of $α_0$ and $β_0$ is $m_1 → m_i → m_i∕(α_iβ_i^{-1})$.
\end{proof*}

We are now ready to give our description of the atoms of $\Sh(\cC^\op,J_\at)$. We will actually describe the quotients of the representable atoms, but since every atom is the image of a representable one, they are all of this form.

\begin{prop}{}{all-atoms}
	Suppose that $\cC$ satisfies \ref{cond:1}, \ref{cond:2} and \ref{cond:3}. Let $n→a$ be a morphism between atoms of $\Sh(\cC^\op,J_\at)$ where $n$ is representable. Then $n→a$ is equal to some composite
	\[ n→m→m∕G ≅ a \]
	where $m$ is another representable atom and $G ⊆ \Aut(m)$.
\end{prop}

\begin{proof*}{}
	Let $n→m$ be a maximal representable quotient such that $n→a$ factorizes through it. Since $\cC^\op$ is supposed co-well-founded by \ref{cond:3}, such an $m$ exists. Let $G ⊆ \Aut(m)$ be the subgroup of automorphisms fixing $m→a$. We claim that $m∕G = a$. Let $x$ be a representable atom and let $α,β : x⇉m$ be two morphisms coequalized by $m→a$. By Lemma~\ref{lem:coeq}, the coequalizer of $α$ and $β$ is of the form $m→m'→m'∕σ$ with $m'$ a representable atom and $σ ∈ \Aut(m')$. But $m→m'$ must be an isomorphism since $m$ is maximal, so that we can suppose $m = m'$ and the coequalizer is $m→m∕σ$ with $σ ∈ G ⊆ \Aut(m)$. Hence we can factorize $m → \coeq(α,β) → m∕G$. But $a$ is the wide pushout of all these coequalizers $\coeq(α,β)$ for $α$ and $β$ coequalized by $m→a$, so that $a = m∕G$.
\end{proof*}

\begin{thm}{}{lfp}
	If $\cC$ satisfies \ref{cond:1}, \ref{cond:2}, \ref{cond:3} and \ref{cond:4}, then $\Sh(\cC^\op,J_\at)$ is locally finitely presentable.
\end{thm}

\begin{proof*}{}
	We must show that any $ω$-chain of atoms stabilizes. Let $n → a_1 → a_2 → ⋯$ be an $ω$-chain of atoms. Suppose that $n$ is representable without loss of generality. By Proposition~\ref{prop:all-atoms}, we can write $n → a_1$ as $n → m_1 → m_1∕G_1 = a_1$. We iterate this process with $m_1 → a_2$ to obtain $m_1 → m_2 → m_2∕G_2 = a_2$. We obtain the following diagram.
	\[\begin{tikzcd}
		n \ar[r] & m_1 \ar[r] \ar[d] & m_2 \ar[r] \ar[d] & m_3 \ar[r] \ar[d] & ⋯\\
		& m_1∕G_1 \ar[r] & m_2∕G_2 \ar[r] & m_3∕G_3 \ar[r] & ⋯
	\end{tikzcd}\]
	Since $\cC$ is well-founded, the sequence $(m_i)_i$ stabilizes at some $N$. Suppose $m_N = m_{N+i}$ for all $i$, without loss of generality. Then $G_N ⊆ G_{N+1} ⊆ G_{N+2} ⊆ ⋯$ and by \ref{cond:4}, this sequence also stabilizes. Hence the category of atoms of $\Sh(\cC^\op,J_\at)$ is co-well-founded. By Proposition~\ref{prop:stab-lfp}, $\Sh(\cC^\op,J_\at)$ is locally finitely presentable.
\end{proof*}

\subsection{Sheaves and pullback-preserving presheaves}
\label{subsec:pBpresSheaves}

We obtained in Proposition~\ref{prop:all-atoms} a way of building all the atoms of $\Sh(\cC^\op,J_\at)$ under the hypothesis \ref{cond:1}, \ref{cond:2} and \ref{cond:3}. However, we didn't give any description of the morphisms between these atoms. In particular, two quotients $n∕G$ and $m∕H$ could be isomorphic as atoms of $\Sh(\cC^\op,J_\at)$ without $n$ and $m$ being isomorphic in $\cC$. This happens for instance in the Malitz--Gregory topos of §~\ref{sec:exmp-malitz-gregory}, as explained in Remark~\ref{rmq:quotients-auto-weird}. We present here in Theorem~\ref{thm:pb-pres-presheaves} a condition under which the category of atoms of $\Sh(\cC^\op,J_\at)$ is easy to describe, as is the case when $\cC = \cFinSetInj$ for instance.

This section would be greatly simplified if one could use a notion of Fraïssé limit of the category $\cC$. A Fraïssé limit of $\cC$ in the naive sense might not exist, since $\Sh(\cC^\op,J_\at)$ might not have any point, as is the case in §~\ref{sec:exmp-malitz-gregory}. However, it seems possible to define something analogous, at least at a formal level. We indicate the intuitions in terms of the Fraïssé limit whenever it is possible, but we give more elementary proofs and we leave the formalization of this concept for future work.

Suppose that $\cC$ satisfies \ref{cond:1} and \ref{cond:2}. We say that $f : A→B$ in $\cC$ is \emph{definable by self-intersections} if every arrow $u : Y→B$ satisfying the following condition factors through $f$:
\begin{center}
	For all $α,β : B⇉X$ equalized by $f$, there is $v : Y→B$ such that $uβ = vα$.
\end{center}
In terms of the Fraïssé limit $\Fc$ of $\cC$, this means that for every point $p$ of $B$ which is not in $A$, there is an automorphism of $\Fc$ which fixes $A⊆B$ and sends $p$ outside of $B$. In general, we will work with partial automorphisms of $\Fc$ defined on a finitary structure $B⊆\Fc$. We represent them by a pair of arrows $α, β : B⇉X$. The arrow $α$ represents the identity embedding $B⊆\Fc$ and $β$ represents the partial automorphism itself.

The lemma below is used to extend automorphisms of $\Fc$ defined on $A⊆\Fc$ to a larger $B$.

\begin{lem}{}{ext-autom}
	Let $\cC$ be a small category with amalgamation. Each diagram of the form on the left can be completed as on the right with $αf' = fα'$ and $βf' = fβ'$.
	
	\begin{minipage}{0.5\textwidth}\[\begin{tikzcd}
		A \ar[r,"α",shift left=0.2em] \ar[r,"β"',shift left=-0.2em] \ar[d,"f"'] & X\\
		B &
	\end{tikzcd}\]\end{minipage}\begin{minipage}{0.5\textwidth}\[\begin{tikzcd}
		A \ar[r,"α",shift left=0.2em] \ar[r,"β"',shift left=-0.2em] \ar[d,"f"'] & X \ar[d,"f'"]\\
		B \ar[r,"α'",shift left=0.2em] \ar[r,"β'"',shift left=-0.2em] & Y
	\end{tikzcd}\]\end{minipage}
\end{lem}

\begin{proof*}{}
	First use amalgamation to build $Y_α$ as below.
	\[\begin{tikzcd}
		A \ar[r,"α"] \ar[d,"f"'] & X \ar[d,"v"]\\
		B \ar[r,"u"] & Y_α
	\end{tikzcd}\]
	Use amalgamation a second time to build $Y_{αβ}$ as follows.
	\[\begin{tikzcd}
		A \ar[r,"β"] \ar[d,"f"'] & X \ar[r,"v"] & Y_u \ar[d,"v'"]\\
		B \ar[rr,"u'"] & & Y_{αβ}
	\end{tikzcd}\]
	The diagram below gives the desired result.
	\[\begin{tikzcd}
		A \ar[r,"α",shift left=0.2em] \ar[r,"β"',shift left=-0.2em] \ar[dd,"f"'] & X \ar[d,"v"]\\[-1.5em]
		& Y_α \ar[d,"v'"]\\[-1.5em]
		B \ar[r,"uv'",shift left=0.2em] \ar[r,"u'"',shift left=-0.2em] & Y_{αβ}
	\end{tikzcd}\]
	Indeed, we have $α(vv') = f(uv')$ and $β(vv') = fu'$.
\end{proof*}

\begin{thm}{}{pb-pres-presheaves}
	Let $\cC$ be a category satisfying \ref{cond:1}, \ref{cond:2} and \ref{cond:3}. Then the following are equivalent:
	\begin{enumerate}
		✦ Every morphism in $\cC$ is definable by self-intersections. \label{itm:pb-def}
		✦ For every representable atom $n$ of $\Sh(\cC^\op,J_\at)$ and every subgroup $G ⊆ \Aut(n)$, the quotient $n∕G$ computed in $[\cC,\cSet]$ is a sheaf. \label{itm:quotient-presh}
		✦ The sheaves for the atomic topology are exactly the pullback-preserving presheaves $\cC→\cSet$. \label{itm:sheaves-pb-pres}
	\end{enumerate}
\end{thm}

Before giving the proof, we must point out that the computations for the implication \ref{itm:quotient-presh}$⇒$\ref{itm:pb-def} were obtained by translating a much clearer argument using Fraïssé limits and presented in Remark~\ref{rmq:real-args-fraisse}. This part of the proof should only be considered as a verification that this translation is possible. Hopefully, such a hack can be eliminated with better tools.

\begin{proof*}{}
	\proofstep{\ref{itm:pb-def}$⇒$\ref{itm:quotient-presh}} Let $A ∈ \cC$ and let $G ⊆ \Aut(B)$. We wish to show that $\Hom(A,—)∕G$ is a sheaf. Let $q : S→T$ in $\cC$, we show that $\Hom(A,—)∕G$ satisfies the sheaf condition relatively to $q$. First, $\Hom(A,S)∕G → \Hom(A,T)∕G$ is injective because $\Hom(A,S)→\Hom(A,T)$ is an injective morphism of $G$-sets. Let $f : A→T$ such that for all pair of morphisms $α,β : T⇉X$ verifying $qα = qβ$, there is some $σ ∈ G$ with $fα = σfβ$. Then $f$ factors through $q$ because $q$ is definable by self-intersections. This shows that $\Hom(A,—)∕G$ is a sheaf.
	
	\proofstep{\ref{itm:quotient-presh}$⇒$\ref{itm:pb-def}} Let $f : A→B$ be a morphism in $\cC$. Our first task is to build a sub-object $j : K→B$ such that that a morphism $p : Y→B$ factors through $K$ exactly when for all $α,β : B⇉X$ equalized by $f$, there is some $p' : Y→B$ such that $pβ = p'α$. 
	We will define inductively a finite sequence of morphisms $j_n : K_n → B$. We start with $j_0 : K_0→B$ equal to the identity of $B$. As long as there are two arrows $α,β : B⇉X$ equalized by $f$ with no $j_n' : K_n→B$ such that $j_nβ = j_n'α$, we form the pullback
	\[\begin{tikzcd}
		K_n \ar[r,"j_n"] & B \ar[r,"β"] & X \\
		K_{n+1} \ar[u] \ar[rr] & & B \ar[u,"α"']
	\end{tikzcd}\]
	and we define $j_{n+1}$ as the composite $K_{n+1}→K_n→B$. The morphism $K_{n+1}→K_n$ is not an isomorphism, otherwise the hypothesis on $α,β$ would be contradicted. Since $\cC$ is well-founded, the sequence of $K_n$ stabilizes at some step $N$ and we define $K = K_N$. We can now show that $j : K→B$ satisfies the desired property. By assumption, for every $α,β : B⇉X$ equalized by $f$, there is some $j' : K→B$ such that $jα=j'β$. Let $p : Y→B$ be any morphism such that for all $α,β : B⇉X$ equalized by $f$, there is some $p' : Y→B$ such that $pβ = p'α$. By induction, $p$ factors through every $K_n$, hence $p$ factors through $K$. We obtain a decomposition of $f$ as
	\[\begin{tikzcd}
		A \ar[r,"i"] & K \ar[r,"j"] & B\text{.}
	\end{tikzcd}\]
	Let $G ⊆ \Aut(K)$ be the set of automorphisms $σ : K→K$ such that $iσ = i$. We will show that the canonical morphism $\Hom(K,—)∕G → \Hom(A,—)$ is an isomorphism. By hypothesis, $\Hom(K,—)∕G$ is a sheaf; in more generally, it is the sheafification of $\Hom(K,—)∕G$ that is isomorphic to $\Hom(A,—)$. The morphism $\Hom(K,—)∕G → \Hom(A,—)$ is an epimorphism in $\Sh(\cC^\op,J_\at)$ since it is a morphism between atoms. It remains to show that $\Hom(K,X)∕G → \Hom(A,X)$ is injective for all $X ∈ \cC$.
	
	Let $α_1,β_1 : K⇉X_1$ such that $iα_1 = iβ_1$. We must show that there is $σ ∈ G$ such that $α_1σ = β_1$. Use Lemma~\ref{lem:ext-autom} to build $X_2$ as below.
	\[\begin{tikzcd}
		A \ar[r,"i"] & K \ar[r,"j"] \ar[d,"α_1"',shift right=0.2em] \ar[d,"β_1",shift left=0.2em] & B \ar[d,"α_2"',shift right=0.2em] \ar[d,"β_2",shift left=0.2em]\\
		& X_1 \ar[r,"j'"] & X_2
	\end{tikzcd}\]
	Then $ijα_2 = ijβ_2$, so that there is $v : K→B$ such that $jβ_2 = vα_2$. We will now show that $v : K→B$ factors through $j : K→B$. That is, we need to show that for every $α_3,β_3 : B→X_3$ equalizing $ij$, there is some $v' : K→B$ such that $vβ_3 = v'α_3$. Use Lemma~\ref{lem:ext-autom} once more to build $X_4$ below.
	\[\begin{tikzcd}
		B \ar[r,"α_3"',shift right=0.2em] \ar[r,"β_3",shift left=0.2em] \ar[d,"α_2"',shift right=0.2em] & X_3 \ar[d,"j''"] \\
		X_2 \ar[r,"α_4"',shift right=0.2em] \ar[r,"β_4",shift left=0.2em] & X_4
	\end{tikzcd}\]
	We have $ijα_2α_4 = ijβ_2β_4$, hence by the universal property of $j:K→B$, there is $v' : K→B$ with $jβ_2β_4 = v'α_2α_4$. We obtain
	\[ vβ_3j'' = vα_2β_4 = jβ_2β_4 = v'α_2α_4 = v'α_3j''\text{.} \]
	Since $j''$ is a monomorphism, we have $vβ_3 = v'α_3$ as desired.
	
	We conclude that there is $σ : k→k$ such that $σj = v$.
	\[\begin{tikzcd}
		K \ar[d,"σ"'] \ar[rd,"v"] & \\
		K \ar[r,"j"] \ar[d,"α_1"',shift right=0.2em] \ar[d,"β_1",shift left=0.2em] & B \ar[d,"α_2"',shift right=0.2em] \ar[d,"β_2",shift left=0.2em] \\
		X_1 \ar[r,"j'"] & X_2
	\end{tikzcd}\]
	Now, $σα_1j' = σjα_2 = vα_2 = jβ_2 = β_1j'$, and since $j'$ is a monomorphism, we have $σα_1 = β_1$. We can see that $σ$ is an automorphism in two ways:
	\begin{enumerate}
		✦ Every endomorphism is an automorphism because $\cC$ is well-founded.
		✦ By symmetry, there is some $σ' : K→K$ such that $σ'β_1 = α_1$. Since $α_1, β_1$ are sub-objects, $σ$ and $σ'$ are inverse isomorphisms.
	\end{enumerate}
	Moreover, $iσα_1 = iβ_1 = iα_1$, so that $iσ = i$. This shows that $σ ∈ G$ and concludes the proof that $\Hom(K,—)∕G→\Hom(A,—)$ is an isomorphism.
	
	Given $u ∈ \Hom(K,X)$, let $[u]$ be its equivalence class modulo $G$. À la Yoneda, let $[i']$ be the image of $\id_A ∈ \Hom(A,A)$ by the inverse morphism $\Hom(A,—)→\Hom(K,—)∕G$, with $i' : K→A$. Then $ii' = \id_A$ and $[i'i] = [\id_K]$ in $\Hom(K,K)∕G$. This means that there is $σ ∈ G$ such that $i'i = σ$. But since $i$ is stabilized by $G$, we have $i'i = i'iσ^{-1} = \id_K$. In other words, $i$ is an isomorphism, which shows that $f : A→B$ is definable by self-intersections.
	
	
	\proofstep{\ref{itm:quotient-presh}$⇒$\ref{itm:sheaves-pb-pres}} We will show that any pullback-preserving presheaf $F : \cC→\cSet$ can be written as a sum of presheaves of the form $n∕G$ with $n$ a representable presheaf and $G ⊆ \Aut(n)$. Note that if $F$ preserves pullbacks, then it sends every morphism to an injection because every morphism of $\cC$ is monic. Given $p ∈ F(X)$, the \emph{support} of $p$ is the minimal $Y ⊆ X$ such that $p$ is in $F(Y) ⊆ F(X)$. Since $\cC$ is well-founded, such a minimal $Y$ exists, and since $F$ preserves pullbacks it is unique. Let $\Ac_F$ be the set of pairs $(p,X)$ with $p ∈ F(X)$ such that the support of $p$ is $X$. Let $G_p$ be the group of automorphisms $σ : X→X$ fixing $p$. For each $(p,X) ∈ \Ac_F$, the element $p$ can be seen as a morphism $p : \Hom(X,—)→F$ and it factors through $\Hom(X,—)∕G_p→F$. We claim that this gives a decomposition
	\[ F ≅ ∑_{(p,X) ∈ \Ac_F} \Hom(X,—)∕G_p \text{.} \]
	Let $Y ∈ \cC$ and let $q : \Hom(Y,—)→F$. We must show that it factors uniquely through one of the $\Hom(X,—)∕G_p→F$. First, the support of $q$ is unique, hence there is a unique $p:\Hom(X,—)→F$ in $\Ac_F$ through which $q$ factors. Since $\Hom(Y,—)$ is projective, $q$ can only factor through one $\Hom(X,—)∕G_p→F$. It remains to show that the factorization of $q$ through $p$ is unique modulo $G_p$. Let $f,g : X⇉Y$ such that $F(f)(p) = F(g)(p)$. Take the following pullback.
	\[\begin{tikzcd}
		X' \ar[r,"u"] \ar[d,"v"] & X \ar[d,"f"]\\
		X \ar[r,"g"] & Y
	\end{tikzcd}\]
	Then $u$ and $v$ are isomorphisms, otherwise the support of $p ∈ F(X)$ would be smaller than $X$. We get that $σ = v^{-1}u$ is an automorphism fixing $p$ and such that $σf = g$. This shows that the factorization of $q$ through $p$ is unique and concludes the proof of the decomposition.
	
	\proofstep{\ref{itm:sheaves-pb-pres}$⇒$\ref{itm:quotient-presh}} We will show that given any pullback-preserving presheaf $F : \cC→\cSet$ and any automorphism group $G ⊆ \Aut(F)$, the quotient $F∕G$ in $[\cC,\cSet]$ is also a pullback-preserving presheaf. First, $F$ takes values in injections because every arrow in $\cC$ is a monomorphism. Any $f : A→B$ in $\cC$ is a monomorphism, hence $F(f) : F(A)→F(B)$ is an injective morphism of $G$-sets and $F(f)∕G : F(A)∕G → F(B)∕G$ is injective too. Consider a pullback in $\cC$.
	\[\begin{tikzcd}
		A \ar[r,"f"] \ar[d,"g"] & B \ar[d,"u"] \\
		C \ar[r,"v"] & D
	\end{tikzcd}\]
	The morphism $F(A)∕G → F(B)∕G ×_{F(D)∕G} F(C)∕G$ is injective since $F(A)∕G → F(B)∕G$ is injective. To show that it is surjective, let $b ∈ F(B)$ and $c ∈ F(C)$ such that there is $σ ∈ G$ with $σ·F(u)(b) = F(v)(c)$. Then $F(u)(σ·b) = F(v)(c)$ so there is $a ∈ F(A)$ with $F(f)(a) = σ·b$ and $F(g)(a) = c$. This shows that $F(A)∕G → F(B)∕G ×_{F(D)∕G} F(C)∕G$ is surjective. As a consequence, for any atom $n$ of $\Sh(\cC^\op,J_\at)$ and any $G ⊆ \Aut(n)$, the quotient $n∕G$ computed in $[\cC,\cSet]$ preserves pushouts, hence is a sheaf.
\end{proof*}

\begin{rmq}{}{real-args-fraisse}
	In terms of Fraïssé limits, the proof of \ref{itm:quotient-presh}$⇒$\ref{itm:pb-def} works as follows. First, we take $A⊆B⊆\Fc$. We define $K⊆B$ as $\setst{p∈B}{∀σ∈\Aut(\Fc) : A⊆\fix(σ) ⇒ σ(p) ∈ B}$. Then any $σ ∈ \Aut(\Fc)$ fixing $A$ satisfies $σ[K] = K$. Let $G ⊆ \Aut(K)$ be the group of automorphisms fixing $A⊆K$. We wish to show that $\Hom(K,X)∕G→\Hom(A,X)$ is injective for all $X$. Let $α,β : K⇉X$ be a pair of morphisms coinciding on $A⊆K$. We extend the embedding $K→\Fc$ to $X$ via $α$, so that $α$ becomes the identity. We choose an automorphism $σ ∈ \Aut(\Fc)$ such that $σ(p) = β(p)$ for all $p ∈ K$. Then $σ$ fixes $A$, so that it restricts to an automorphism $σ' : K→K$. Thus $β(p) = σ(p) = σ'(p)$ for all $p ∈ K$. Recall that $α$ was identified with the identity injection, so we get $β = α∘σ'$. This concludes the proof.
	
	In a similar way, we can show that the morphisms definable by self-intersections are stable under composition. Using Fraïssé limits, the argument is simple: Let $A⊆B⊆C⊆\Fc$ where $A⊆B$ and $B⊆C$ are definable by self-intersections. Let $p ∈ C⧵A$. If $p ∉ B$, then there is an automorphism of $\Fc$ fixing $B$, hence $A$, and sending $p$ outside of $C$. Otherwise, $p ∈ B⧵A$ and there is $σ ∈ \Aut(\Fc)$ fixing $A$ and sending $p$ outside of $B$. If $σ(p) ∉ C$, we are done. Otherwise, $σ(p) ∈ C⧵B$ and there is another automorphism $σ'$ of $\Fc$ fixing $B$ and such that $σ'(σ(p)) ∉ C$. In this case, $σ'∘σ$ fixes $n$ and sends $p$ outside of $C$.
	
	This argument can also be transformed into a proof ``without Fraïssé limits,'' at the cost of obscuring the general idea.
\end{rmq}

\begin{rmq}{}{}
	Every morphism of $\cC$ definable by self-intersections is the equalizer of all the pairs of morphisms that it equalizes. This means that half of \ref{cond:1} is contained in \ref{itm:pb-def} of Theorem~\ref{thm:pb-pres-presheaves}. In particular, $f : A→B$ is definable by self-intersections when there is a pair of morphisms $α,β : B⇉X$ such that the square below is cartesian.
	\[\begin{tikzcd}
		A \ar[r,"f"] \ar[d,"f"'] & B \ar[d,"α"]\\
		B \ar[r,"β"] & X
	\end{tikzcd}\]
\end{rmq}

When the conditions of Theorem~\ref{thm:pb-pres-presheaves} are satisfied, the category of atoms is very simple to describe, thanks to \ref{itm:quotient-presh}. Every atom is of the form $n∕G = \Hom(n,—)∕G$ with $n ∈ \cC$, and the morphisms $n∕G→m∕H$ can be described as the arrows $f : m→n$ in $\cC$ such that for every $γ ∈ H$, there is some $σ ∈ G$ with $γf = fσ$. Moreover, two such arrows $f, g : m→n$ give the same morphism $n∕G→m∕H$ if and only if there is $σ ∈ G$ with $fσ = g$.

\begin{quest}{}{}
	This explicit description of $\Sh(\cC^\op,J_\at)$ can be reproduced for any category $\cC$, even when it doesn't describe $\Sh(\cC^\op,J_\at)$. Is there a direct proof that this category is a topos? What hypothesis on $\cC$ are needed?
\end{quest}


\section{The Malitz--Gregory atomic topos}
\label{sec:exmp-malitz-gregory}

The Malitz--Gregory topos is an example of a non-degenerate atomic topos which doesn't have any points. It was first given in \cite[Sec.~5]{makkaiFullContinuousEmbeddings1982}, based on an earlier example in infinitary logic of \cite{HanfNumberCompleteMal1968} and \cite{IncompletenessFormalSystemGre1971}. A description can also be found in \cite[D3.4.14]{johnstoneSketchesElephantTopos2002}.

We define in §~\ref{subsec:siteMG} an atomic site of presentation of this topos adapted to our study. This site was obtained by first considering a simpler atomic site and by sheafifying the representable presheaves. We show in §~\ref{subsec:MGwellfounded} that this topos is locally finitely presentable by applying the criterion of §~\ref{sec:general-crit}.

\subsection{Definition of the Malitz--Gregory topos}
\label{subsec:siteMG}

A \emph{full binary tree} is a tree in which each node is either a leaf, or has exactly two children. More formally, it can be defined as a pointed oriented graph $(X,r_X ∈ X,{\prec} ⊆ X^2)$ such that:
\begin{enumerate}
	✦ For each $x ∈ X$, the set $\setst{y∈X}{x \prec y}$ is either empty or of cardinality $2$.
	✦ For each node $y$, there is exactly one path $r_X = x_0 \prec x_1 \prec ⋯ \prec x_n = y$.
\end{enumerate}
When $x \prec y$, we say that $y$ is a \emph{child} of $x$ and that $x$ is the \emph{parent} of $y$. An \emph{embedding} of trees $f : X→Y$ is an (injective) graph embedding preserving the root.

\paragraph{$\II$-trees} We now enrich our notion of tree with a partial labeling of its infinite branches. Let $\II$ be set. A \emph{branch} of $X$ is an infinite sequence $r_X = x_0 \prec x_1 \prec x_2 \prec ⋯$ Let $B_X$ be the set of branches of $X$. An \emph{$\II$-tree} is a full binary tree $X$ equipped with a partial function $c_X : B_X → \II$. An embedding of $\II$-trees is an embedding of full binary trees that preserves the labeling function.

An $\II$-tree is \emph{finitary} if it has a finite number of branches, a finite number of elements whose parent is not on a branch, and a total labeling function. Let $\cT_\II$ be the category of finitary $\II$-trees and embeddings. A sub-tree of an $\II$-tree is a sub-object in $\cT_\II$.

\begin{lem}{}{trees-ind}
	The category of $\II$-trees is equivalent to $\cInd(\cT_\II)$.
\end{lem}


\begin{proof*}{}
	First, the category of $\II$-trees has filtered colimits, by taking the union of the underlying sets and labeling functions. Moreover, each $\II$-tree is the filtered union of its finitary sub-trees. We can also show that each finitary $\II$-tree is finitely presentable: Let $⋃_{i∈I} X_i$ be a filtered union of $\II$-trees and let $X ⊆ ⋃_{i∈I} X_i$ be a finitary sub-tree. Each branch of $X$ is labeled, thus it is labeled as a branch of $⋃_{i∈I} X_i$. Since the labeling function of $⋃_{i∈I} X_i$ is the union of the labeling functions of the $X_i$, the branch must be contained in one of the $X_i$. For each element $x$ of the branch, one of its children is contained in $X_i$, hence the other one too. Only a finite number of elements of $X$ are not covered in this way, and since the union $⋃_{i∈I} X_i$ is filtered, $X$ is a sub-tree of one of the $X_i$.
\end{proof*}

We impose the additional condition that the cardinality of $\II$ is strictly greater than $2^{ℵ_0}$. The \emph{Malitz--Gregory topos} over $\II$ is $\Sh(\cT_\II^\op,J_\at)$.

\begin{prop}{}{pointless}
	$\Sh(\cT_\II^\op,J_\at)$ has no points.
\end{prop}

\begin{proof*}{}
	By Diaconescu's theorem, a model of $\Sh(\cT_\II^\op,J_\at)$ in $\cSet$ is a flat functor $\cT_\II^\op → \cSet$ continuous with respect to the atomic topology, i.e., which sends every morphism to a surjection. Since flat functors $\cT_\II^\op → \cSet$ can be identified with ind-objects of $\cT_\II$, and thanks to Lemma~\ref{lem:trees-ind}, a model is a special $\II$-tree $X ∈ \cInd(\cT_\II)$. The condition that the corresponding functor sends every morphism to a surjection translates into the fact that each morphism $n → X$ with $n ∈ \cT_\II$ extends along any morphism $n→m$ in $\cT_\II$.
	\begin{equation}\label{eq:extension-model}\begin{tikzcd}
			n \ar[r] \ar[d] & X\\
			m \ar[ru,dashed,"∃"'] &
	\end{tikzcd}\end{equation}
	However, no $\II$-tree can satisfy this extension condition. Suppose that $T$ is such a tree. Like any tree as we defined them, $T$ has at most a countable number of nodes. Moreover, for each $i ∈ \II$, there must be a branch labeled by $i$, by applying the extension property \eqref{eq:extension-model} to the inclusion of the one-point tree into an infinite branch labeled by $i$, together with the other nodes forced to be present in order to have a full binary tree. In particular, the partial function $B_T → \II$ must be surjective, but it is not possible because $B_T ⊆ 𝒫(T)$ and $\II$ has cardinality strictly bigger than that of $𝒫(T)$. In conclusion, $\Sh(\cT_\II^\op,J_\at)$ has no model in $\cSet$.
\end{proof*}

\subsection{The Malitz--Gregory topos is locally finitely presentable}
\label{subsec:MGwellfounded}

Given two nodes $x$ and $y$ in a tree $T$, we say that $y$ is a descendant of $x$ or that $x$ is an ancestor of $y$ if there is a path $x = x_0 \prec x_1 \prec ⋯ \prec x_n = y$. Given a node $x ∈ T$, we denote by ${↓_T}x$ the set of all the descendants of $x$ equipped with the induced tree structure. We use simply ${↓}x$ when $T$ is implicit. (It is not a sub-tree because the injection ${↓}x→T$ doesn't preserve the root.)

\begin{lem}{}{subcan}
	The category $\cT_\II$ satisfies \ref{cond:1}.
\end{lem}

\begin{proof*}{}
	First, we show that $\cT_\II$ has amalgamation. This means that given two finitary $\II$-trees $A$ and $B$, and given a common sub-tree $X ⊆ A$, $X ⊆ B$, we can find a tree $C$ containing both $A$ and $B$ as sub-trees, such that the intersection contains $X$. Let $W$ be the complete binary tree of infinite depth. Choose an arbitrary embedding $X ⊆ W$. Extend it arbitrarily to two embeddings $A ⊆ W$ and $B ⊆ W$. We now wish to label some of the branches of $W$ so that $A ⊆ W$ and $B ⊆ W$ become an embeddings of $\II$-trees. This might be impossible if a branch of $A$ and a branch of $B$ with different labels are sent to the same branch of $W$. Nonetheless, we can correct these ``conflicts'' as follows. For each such conflict, there is a branch of $A$ and a branch of $B$ with different labels that get identified in $W$. This path is not in $X$, since the labels would coincide otherwise. Pick a node $x$ of this path which is deep enough to ensure that $x ∉ X$ and that both ${↓_A}x$ and ${↓_B}x$ are composed of a unique branch containing every node or its parent. Modify the embedding of $A$ and $B$ in $W$ from $x$ onward so that the two paths are not sent to the same one anymore. Doing so for each conflict creates two new embeddings $A,B ⊆ W$. The union $A∪B$ is finitary and can be equipped with two labelings so that the embeddings of $A$ and $B$ preserve the labels. This shows that $\cT_\II$ has amalgamation.
	
	We will now show that every morphism of $\cT_\II$ is a regular monomorphism.
	
	Let us show that for each finitary $\II$-tree $A$, there are two embeddings $A ⇉ A'$ whose equalizer contains only the root of $A$. If $A$ contains only a root, we can take $A' = A$. If not, let $x_1$ and $x_2$ be the two children of the root of $A$. Let $C$ be a finitary tree such that ${↓}x_1 ⊆ C$ and ${↓}x_2 ⊆ C$ (this is possible because $\cT_\II$ has amalgamation). Let $A'$ be the tree obtained by joining two copies of $C$ by a root. Then there are two obvious embeddings $A ⇉ A'$ and their equalizer contains only the root of $A$.
	
	Now, let $X ⊆ Y$ be an arbitrary embedding of finitary $\II$-trees. We show it is a regular monomorphism. Let $L = \setst{x∈X}{{↓_X}x=\set{x}}$. Let $Y'$ be the $\II$-tree obtained by replacing each ${↓_Y}x ⊆ Y$ for $x ∈ L$ by the tree $[{↓_Y}x]'$ obtained in the previous paragraph, with two embeddings ${↓}x ⇉ [{↓}x]'$ whose equalizer is the identity. This gives two embeddings $Y⇉Y'$ whose equalizer is $X$, since every element of $Y⧵X$ is a descendant of some element of $L$. Hence $X ↪ Y$ is a regular monomorphism and the atomic topology on $\cT_\II^\op$ is sub-canonical.
\end{proof*}

\begin{rmq}{}{}
	Even though $\cT_\II$ has amalgamation, it doesn't have \emph{functorial} amalgamation as defined in \cite[Dfn.~4.1.6]{dilibertiTopoiEnoughPoints2024}, since it would imply that $\Sh(\cT_\II,J_\at)$ has enough points by \cite[Cor.~4.1.8]{dilibertiTopoiEnoughPoints2024}.
\end{rmq}

\begin{lem}{}{pres-push}
	The category $\cT_\II$ satisfies \ref{cond:2}.
\end{lem}

\begin{proof*}{}
	We show that $\cT_\II$ satisfies the equivalent condition \ref{cond:2-}. Consider a pullback as below in $\cT_\II$.
	\[\begin{tikzcd}
		X∩Y \ar[r,hook] \ar[d,hook] & X \ar[d,hook] \\
		Y \ar[r,hook] & Z
	\end{tikzcd}\]
	Let $u, v : Z⇉A$ be two parallel arrows coinciding on $X∩Y$. Let $L = \setst{x∈X∩Y}{{↓_X}x=\set{x}}$. Note that an element of $L$ cannot be a descendant of another element of $L$. Define $w : L→A$ by
	\[w(x) = \begin{cases*}
		v(x) & if $x$ has an ancestor in $L$,\\
		u(x) & if $x$ has no ancestor in $L$.
	\end{cases*}\]
	Then $w$ is a tree embedding because it is obtained by modifying the definition of $u$ on ${↓}_Zx$ for each $x∈L$ to fit another tree embedding $v$ having the same definition on $x$. Moreover:
	\begin{enumerate}
		✦ $w$ and $u$ coincide on $X$ because if $x∈X$ has an ancestor in $L$, then this ancestor is $x$ itself and $u(x) = v(x)$;
		✦ $w$ and $v$ coincide on $Y$ because for each $y∈Y$, either $y∈X∩Y$ and $u(y)=v(y)$, or $y∉X∩Y$ and the closest ancestor of $y$ in $X∩Y$ is in $L$.
	\end{enumerate}
	This shows that $\cT_\II$ satisfies \ref{cond:2-} hence \ref{cond:2}.
\end{proof*}

\begin{lem}{}{finiteness}
	$\cT_\II$ satisfies \ref{cond:3} and \ref{cond:4}.
\end{lem}

\begin{proof*}{}
	We must show that for every finitary $\II$-tree $X$, the poset of sub-trees of $X$ is well-founded. Recall that $B_X$ is the set of branches of $X$. Let $F_X$ be the set of elements of $X$ whose parent is not in a branch of $X$. If $Y ⊆ X$ is a proper sub-tree, then $\card{B_Y} ≤ \card{B_X}$, and if $\card{B_Y} = \card{B_X}$, then $\card{F_X} < \card{F_X}$. Hence $X ↦ (\card{B_X}, \card{F_X})$ defines a strictly order-preserving map $\cT_\II → ω^2$ and $\cT_\II$ is well-founded.
	
	To show that the groups of automorphisms of $\cT_\II$ are Noetherian, note that any morphism $X→Y$ in $\cT_\II$ is uniquely determined by a function $B_X → B_Y$ and a function $F_X → F_Y$, thus there are only finitely many of them. The automorphism groups are thus finite, and in particular Noetherian.
\end{proof*}

Using Theorem~\ref{thm:lfp}, we obtain:

\begin{cor}{}{}
	$\Sh(\cT_\II^\op,J_\at)$ is a connected locally finitely presentable atomic topos with no points.
\end{cor}

\begin{rmq}{}{quotients-auto-weird}
	We know by Proposition~\ref{prop:all-atoms} that every atom of $\Sh(\cT_\II^\op,J_\at)$ can be described as a formal quotient of a finitary $\II$-tree $n$ by a group of automorphisms $G$. However, contrary to what happens with the Schanuel topos, it is possible that $n∕G$ and $m∕H$ are isomorphic but that $n$ and $m$ are not. For instance, let $n$ be the unique $\II$-tree with exactly $3$ nodes, or rather the atom corresponding to this tree. Let $G$ be the automorphism group of $n$, containing the identity and an involution. Then $n∕G$ is isomorphic to the atom associated to the tree with only one node.
\end{rmq}

\section{Final remarks}
\label{sec:final-rmk}

Lastly, we indicate a list of possible future directions:
\begin{enumerate}
	✦ Describe more explicitly the category of atoms of $\Sh(\cC^\op,J_\at)$ than what we obtained here, as suggested by Remark~\ref{rmq:quotients-auto-weird}.
	✦ Extend the results of this paper to the \emph{prime-generated toposes} of \cite{makkaiFullContinuousEmbeddings1982}.
	✦ Another direction of exploration is the relation with Fraïssé limits and Galois theory. See \cite{FraisseConstructionToposTheoreticCar2014,AspectsTopologicalGaloisCarLaf2019} for more details. The topos $\Sh(\cC^\op,J_\at)$ might have no points, but the intuition of working with the Fraïssé limit of $\cC$ still applies, as illustrated in §~\ref{subsec:pBpresSheaves}. When $\cC$ has a Fraïssé limit, the open subgroups of its automorphism group correspond to the objects of $\cC$. The well-foundedness condition on $\cC$ is thus related to a Noetherianity condition of this automorphism group. 
	✦ Some of the ideas in §~\ref{subsec:pBpresSheaves}, in particular the notion of support, are very similar to the ideas of \cite{NotesCombinatorialFunctorsFio2001,ReflectiveKleisliSubcategoriesFioMen2005} and \cite{PolyadicSetsHomomorphismReg2022}, where well-foundedness also plays an important role. In \cite{PolyadicSetsHomomorphismReg2022}, the author considers \emph{polyadic sets} and their \emph{Stirling kernels} from a combinatorial point of view. A logical interpretation of the same kind of phenomenon is possible, see for instance \cite{PolyadicSpacesProfiniteMar2021} or \cite[Chap.~5]{LogiqueCategoriquePointMar2023}. The logical interpretation of polyadic sets is that they are dual to the geometric hyperdoctrines whose classifying toposes are atomic. Curiously, the pullback-preserving functors $\cC→\cSet$ in Theorem~\ref{thm:pb-pres-presheaves} are also exactly the polyadic sets over $\cC^\op$ (and thus there is an atomic topos naturally associated to each one of them).
\end{enumerate}

\AtNextBibliography{\small}
{\printbibliography[
heading=bibintoc,
title={References}
]}

\end{document}